\documentclass[11pt]{article}
\begin{document}
\tolerance=10000
\newcommand{\eproof}{\rule{0.2cm}{0.2cm}}
\newcommand{\FTS}[2]{\frac{{\textstyle #1}}{{\textstyle #2}}}
\newcommand{\NN}{I\!\!N}
\newcommand{\ka}{I\!\!K}
\newcommand{\rgr}{{\rm grad}}
\newcommand{\CC}{I\!\!\!\!C}
\newcommand{\RR}{I\!\!R}
 \newcommand{\intl}{\int\limits}
\newcommand{\suml}{\sum\limits}
\setcounter{page}{1}
\thispagestyle{empty}
\font\bfs=cmbx10 scaled \magstep2
\def\eg{{\it e.g.}\ } \def\ie{{\it i.e.}\ }
\def\sg{\hbox{sign}\,}
\def\sgn{\hbox{sign}\,}
\def\sign{\hbox{sign}\,}
\def\e{\hbox{e}}
\def\exp{\hbox{exp}}
\def\ds{\displaystyle}
\def\dis{\displaystyle}
\def\q{\quad}    \def\qq{\qquad}
\def\lan{\langle}\def\ran{\rangle}
\def\l{\left} \def\r{\right}
\def\lra{\Longleftrightarrow}
\def\arg{\hbox{\rm arg}}
\def\d{\partial}
 \def\dr{\partial r}  \def\dt{\partial t}
\def\dx{\partial x}   \def\dy{\partial y}  \def\dz{\partial z}
\def\rec#1{{1\over{#1}}}
\def\log{\hbox{\rm log}\,}
\def\erf{\hbox{\rm erf}\,}     \def\erfc{\hbox{\rm erfc}\,}
\def\G{{G_{\alpha,\beta}^\theta}}
\def\K{K_{\alpha,\beta}^\theta}
\def\Gxt{\G (x,t)}
\def\Gkt{{\widehat{\G}}  (\kappa,t)}
\def\Gxs{{\widetilde{\G}}  (x,s)}
\def\Gks{{\widehat{\widetilde {\G}}} (\kappa,s)}
\def\Gos{{\widehat{\widetilde {\G}}} (0,s)}
\def\Got{{\widehat{\G}} (0,t)}
\def\Gone{G_{\beta}^{(1)}} \def\Gtwo{G_{\beta}^{(2)}} \def\Gj{G_{\beta}^{(j)}}
\def\Kone{K_{\beta}^{(1)}} \def\Ktwo{K_{\beta}^{(2)}} \def\Kj{K_{\beta}^{(j)}}
\def\Gxtone{\Gone (x,t)}  \def\Gxttwo{\Gtwo (x,t)} \def\Gxtj{\Gj (x,t)}
\def\Gktone{{\widehat{\Gone}}(\kappa,t)}
\def\Gkttwo{{\widehat{\Gtwo}}(\kappa,t)}
\def\Gktj{{\widehat{\Gj}}(\kappa,t)}
\def\Gxsone{{\widetilde{\Gone}}  (x,s)}
\def\Gxstwo{{\widetilde{\Gtwo}}  (x,s)}
\def\FT{{\cal F}\,} 
\def\LT{{\cal L}\,}  
\def\L{{\cal L}} 
\def\F{{\cal F}} 
\def\M{{\cal M}}  
\def\I{{\cal I}}  
\def\pni{\par \noindent}
\def\vsh{\smallskip}
\def\vs{\medskip}
\def\vvs{\bigskip}
\def\vvvs{\bigskip\medskip} 
\def\vsp{\par}
\def\vsn{\vsh\pni}
\def\cen{\centerline}
\def\ra{\item{a)\ }} \def\rb{\item{b)\ }}   \def\rc{\item{c)\ }}
\newcommand{\lt}{\left}
\newcommand{\rt}{\right}
\newcommand{\R}{{\mathcal R}}
\newcommand{\Mt}{\,\stackrel{\mathcal M}{\longleftrightarrow}\,}
\newcommand{\Ft}{\,\stackrel{\mathcal F}{\longleftrightarrow}\,}
\newcommand{\Lt}{\,\stackrel{\mathcal L}{\longleftrightarrow}\,}
 \cen{FRACALMO PRE-PRINT   {\bf www.fracalmo.org}}
\vsh
\cen{\bf Fractional Calculus and Applied Analysis,
  Vol. 6 No 4 (2003), pp. 441-459}
\vsh
\cen{An International Journal for Theory and Applications \ ISSN 1311-0454}
\vsh
\cen{{\bf www.diogenes.bg/fcaa/}}
\vs
\hrule
   \vskip 0.50truecm


\font\title=cmbx12 scaled\magstep2
\font\bfs=cmbx12 scaled\magstep1
\font\little=cmr10
\begin{center}

{\title Mellin transform and subodination laws}
\vs

{\title in fractional diffusion processes}

\vvs

 {Francesco MAINARDI}$^{(1)}$, {Gianni PAGNINI}$^{(2)}$
and
{Rudolf GORENFLO}$^{(3)}$


\vs

$\null^{(1)}$
 {\little Department of Physics, University of Bologna, and INFN,} \\
{\little Via Irnerio 46, I-40126 Bologna, Italy} \\
{\little Corresponding Author; E-mail: {\tt mainardi@bo.infn.it}}  
\\ [0.25 truecm]
$\null^{(2)}${\little ENEA: Italian Agency for  New Technologies,
  Energy and the Environment} 
\\ {\little Via Martiri di Monte Sole 4, I-40129 Bologna, Italy}
   \\ [0.20 truecm]
$\null^{(3)}$
 {\little Department of Mathematics and  Informatics,
 Free   University of Berlin,} \\
 {\little Arnimallee  3, D-14195 Berlin, Germany} 

 \def\date#1{\gdef\@date{#1}} \def\@date{\today}

\end{center}

\vskip 0.10truecm
\noindent
{\it This paper is dedicated to Paul Butzer,   Professor Emeritus, 
Rheinisch-Westf\"alische Technische Hochschule  (RWTH), Aachen, Germany,
on the occasion of his 75-th birthday (April 15, 2003)}

\vs

\cen{\bf Abstract} 
\vs
\noindent
The Mellin transform is usually applied
in probability theory    to the product of
independent  random variables.
In recent times  the machinery of the Mellin transform
has been adopted to describe  the L\'evy stable  distributions,
and more generally   the probability distributions
governed by generalized diffusion equations of fractional order
in space and/or in time.
In these cases the related stochastic processes
are self-similar and are  simply referred to as
fractional diffusion processes.
In this note, by using the convolution properties of the Mellin transform,
we  provide some (interesting) integral formulas
involving the distributions of these processes
that can be interpreted  in terms of  {\it subordination laws}.



\section*{1. Introduction}

The role of the Mellin transform in probability theory is
mainly  related  to the product of
independent  random variables:  in fact
it is well-known that the probability density of the product of
two independent  random variables is given by  the Mellin convolution of
the two corresponding densities.
Less known is  their role with respect to the class
of the L\'evy stable distributions, that was formerly outlined
by Zolotarev \cite{Zolotarev 86} and Schneider \cite{Schneider LNP86},
see also \cite{UchaikinZolotarev 99}.
A  general class of probability distributions (evolving in time),
that includes the L\'evy strictly stable distributions,
is obtained by solving, through  the machinery of  the Mellin transform,
 generalized  diffusion equations of fractional order  in space
and/or  in time, see  \cite{GoIsLu 00,Mainardi LUMAPA01}.
In these cases the related stochastic processes
turn out  as self-similar and are
referred to as fractional diffusion
processes.

In this note, after  the essential notions and notations
concerning  the  Mellin transform,
we  first show the role of the Mellin convolution
between  probability densities to establish  subordination laws
related to  self-similar stochastic processes.
Then,  for the fractional diffusion processes
we establish   a sort of Mellin convolutions
between the related  probability densities,
that can be interpreted as  {\it subordination laws}.
This is carried out starting for the
representations through  Mellin-Barnes integrals
of the probability densities.

We point out that our  results, being based on simple manipulations, 
can be  understood by non-specialists of transform methods and
special functions; however they
could be derived through a more general analysis
involving the class of
higher transcendental functions
of Fox $H$ type to which the
probability densities arising as fundamental solutions
 of the fractional diffusion equation  belong.



\section*{2. The Mellin transform}

The Mellin  transform of a sufficiently well-behaved function $f(x)$,
$x \in \RR^+$,
is defined by
$$
   {\cal M} \, \{ f(x); s\} = f^*(s)=
   \int_0^{+\infty} f(x)\,
 x^{s-1}\,  dx\,, \q s\in \CC \,,\ 
\eqno(2.1)
$$
when the integral converges.
Here we  assume $f(x) \in L_{loc}(\RR^+)$ according
to the most usual approach suitable for applied scientists.
The basic properties of the Mellin transform follow immediately
from those of the bilateral Laplace transform
since the two transforms are intimately connected.

Recently the theory of  the  Mellin transform   has been
the object of  intensive researches by Professor  Butzer
and his associates, see \eg
\cite{ButzerJansche 97,ButzerJansche 98,
Butzer MellinFC1-02,Butzer MellinFC2-02,Butzer MellinFC3-02,
ButzerWestphal 00};
in particular
Butzer and Jansche \cite{ButzerJansche 97,ButzerJansche 98}
have introduced a  theory
independent from Laplace or Fourier transforms.

The integral (2.1) defines the Mellin transform in a vertical strip
in the $s$-plane whose boundaries are determined by the analytic structure
of $f(x)$ as $x \to 0^+$ and $x \to +\infty$.

If we suppose that
$$ f(x)= \cases{
O\l(x^{-\gamma _1 -\epsilon }\r)
  & as $\; x \to 0^+\,, $\cr\cr
O\l(x^{-\gamma _2 -\epsilon}\r)
  & as $\; x \to +\infty\,, $\cr}  \eqno(2.2)
$$
for every (small) $\epsilon>0 $
and $\gamma _1<\gamma _2$, the the integral (2.1)
converges absolutely and defines an analytic function
in the strip
$ \gamma _1 < \Re {s} < \gamma _2 \,. $
This strip is known as the {\it strip of analyticity} of
${\cal M} \, \{ f(x ); s\} = f^*(s)$

The inversion formula for (2.1) follows directly from the corresponding
inversion formula for the bilateral Laplace transform. 
We have
$$
{\M}^{-1}\{f^*(s);x\}=f(x)=
\frac{1}{2 \pi i} \int_{\gamma -i \infty}^{\gamma +i \infty}
 f^*(s)\, x^{-s}\, ds \,,\q 
\gamma_1 < \gamma < \gamma_2\,,
\eqno(2.3)$$
at all points $x \ge 0$ where $f(x)$ is continuous.

Let us now consider the most relevant {\it Operational Rules}.
Denoting by   $\stackrel{{\cal M}} {\leftrightarrow}$
the juxtaposition of a function $f(x)$ on $x >0$ with its
Mellin transform $f^*(s)\,,$  
we have
$$
x^a \, f(x) \Mt f^*(s+a) \,, \q a \in \CC\,,
$$
$$
   f(x^b) \Mt \frac{1}{|b|} f^*(s/b)\,, \q b \in \CC\,,\q  b\neq 0\,,
$$
$$
f(cx) \Mt c^{-s}\, f^*(s) \,, \q  c \in \RR\,,\q  c>0 \,,
$$
from which
$$
x^a \, f(c x^b) \Mt \frac{1}{|b|}
\, c^{- (s+a)/b}\, f^*\lt(\frac{s+a}{b}\rt)\,. 
\eqno(2.4)$$
Furthermore we have
$$
h(x)=\int_0^{\infty} f \lt(\frac{x}{\xi} \rt)\, g(\xi)\, \frac{d\xi}{\xi}
\Mt
 f^*(s)\, g^*(s) = h^*(s) \,,
\eqno(2.5)$$
which is known as
the {\it Mellin convolution} formula.


\section*{3. Subordination in stochastic processes via Mellin convolution}

In recent years a number of  papers have appeared
where explicitly or implicitly
subordinated stochastic processes
have been treated in view of their relevance
in physical and financial applications, see e.g.
\cite{BaeumerMeerschaert 01,Barkai ChemPhys02,OEBN EDITOR01,Nevada PRE02,
Metzler-Klafter PhysRep00,Scalas PRE04,Stanislavsky PHYSA03,
UchaikinZolotarev 99,Wyss-Wyss 01}
and references therein.
Historically, the notion of subordination was originated by Bochner,
see \cite{Bochner 55,Bochner 62}.
In brief, according to Feller \cite{Feller 71},
 a {\it subordinated process} $X(t)= Y(T(t))$
is obtained by randomizing the time clock of
a stochastic process $Y(t)$ using a new clock $T(t)$,
where $T(t)$ is a random process with non-negative independent
increments. The resulting process $Y(T(t))$ is said
to be subordinated to $Y(t)$, called
the {\it parent process}, and is directed by $T(t)$ called
the {\it directing process}.
The directing process is often referred to as
the  operational time\footnote{
ADDED NOTE (2007). In a recent paper by R. Gorenflo, F. Mainardi and A. Vivoli:
Continuous time random walk and parametric subordination in fractional
diffusion ({\tt http://arXiv.org/abs/cond-mat/0701126}),
to be published in Chaos, Solitons \& Fractals (2007),
the Authors
have given an interesting alternative  interpretation of subordination, that they call
{\it parametric subordination}.}.
In particular, assuming  $Y(\tau)$ to be a Markov process
with a spatial probability density function ($pdf$)
of $x$, evolving in time $\tau$,
$q_\tau (x) \equiv q(x; \tau)$, and  $T(\tau)$ to be a 
process with non-negative independent increments
with $pdf$ of $\tau$ depending on a parameter $t$,
$u_t(\tau) \equiv u(\tau; t)$, then the subordinated
process $X(t) = Y(T(t))$ is governed
by the spatial $pdf$ of $x$
evolving with $t$,
$p_t (x) \equiv p(x; t),$  given by the integral representation
$$
p_t(x)=\int_0^{\infty} q_{\tau}(x)\, u_t(\tau)\, {d\tau} \,.
\eqno(3.1)$$
If the parent process $Y(\tau)$ is  {\it self-similar}
of the kind that
 its $pdf$ $q_{\tau}(x)$ is such that
$$
q_{\tau}(x)\equiv q(x;\tau) =
\tau^{-\gamma } \, q \lt(\frac{x}{\tau^{\gamma }}\rt)\,,
\q \gamma >0\,,
\eqno(3.2) $$
then Eq. (3.1) reads,
$$
p_t(x)=\int_0^{\infty} q_{\tau}\lt(\frac{x}{\tau^{\gamma }}\rt)
u_t(\tau) \frac{d\tau}{\tau^{\gamma }} \,.
\eqno(3.3)$$
Herewith  we  show how to interpret Eq. (3.3)
 in terms of a special  convolution integral in the
framework of the theory of  the Mellin transform.
Later, in the next sections, we shall show how to use the tools
 of the Mellin-Barnes integral  and Mellin transform  to
treat the subordination for the  class of self-similar
stochastic process, which are governed by fractional
diffusion equations.

Let $X_1$ and $X_2$
be two 
real {\it independent} random variables with $pdf$'s
$p_1(x_1)$ and $p_2(x_2)$ respectively,
with $x_1\in \RR$ and $x_2 \in \RR_0^+$.
The joint probability is
$$
p_*(x_1,x_2)=p_1(x_1)\, p_2(x_2)\,.
\eqno(3.4)$$
Denoting by $X$  the random variable
obtained by the product of $X_1$ and $X_2^\gamma$,
i.e.  $x=x_1\, x_2^\gamma\,,$
and carrying out the transformation
$$
\cases{
x_1  = x/\tau^\gamma\,,\cr 
x_2 =\tau\,,\cr
}
\eqno(3.5)$$
we get the identity
$$
\widetilde {p_*}(x,\tau)\, dx \, d\tau=p_1(x/\tau^\gamma)\, p_2(\tau)\, J\,
dx \, d\tau\,,
\eqno(3.6)$$
where
$$ J = \l\vert \;
{ {\ds{\d x_1\over \d x }}\q {\ds {\d x_1\over \d x_2 }}
\atop
  {\ds {\d x_2\over \d x }}\q {\ds {\d x_2\over \d x_2 }} }
\; \r\vert
= \l\vert
{ {\ds {1\over \tau^\gamma }} \q - {\ds {\gamma x \over \tau^{\gamma+1} }}
\atop
  {\!\! 0 }\qq \q\q {1 } }
\r\vert 
 \eqno(3.7)$$
 is the Jacobian of  the transformation (3.5).
Noting that
$ J = 1/\tau^\gamma $ and integrating  (3.6)
  in $d\tau$ we finally get the $pdf$ of $X$,
$$
p(x)=\int_{\RR^+}p_1\l(\frac{x}{\tau^\gamma}\r)\,
p_2(\tau)\,\frac{d\tau}{\tau^\gamma}\,, \q x\in \RR\,.
\eqno(3.8)$$
For $\gamma=1$, by comparing with Eq. (2.5),
we recover the well known result
that  the probability density of the product of
two  independent  random variables is given by
the Mellin convolution of
the two corresponding densities.

We now adapt Eq. (3.8) to our subordination formula
(3.3) by identifying $p_1, p_2$ with $q_\tau$ and $u_t$,
respectively.
We  can now interpret
the subordination formula (3.3)
as follows.
The $pdf$ of the subordinated process $X$, $p_t(x)$,
turns out to be the $pdf$ of the product of
the independent random  variables $X_q$ and ${X_u}^{\gamma }$
distributed according to $q_\tau(x_q)$ and $u_t(x_u)$,
respectively.


\section*{4. Fractional diffusion equations and probability distributions}

An interesting way to generalize the classical diffusion equation
$${\d^2 \over \d x^2} \,u(x,t) \, = \,
    {\d \over \d t}\, u(x,t)\,, \q -\infty<x<+\infty\,, \q t \ge 0 \,,
\eqno(4.1)  $$
is to replace in (4.1) the partial derivatives in space and time
by suitable linear integro-differential operators,
to be intended as  derivatives of non integer order,
that allow the corresponding  Green function (see below) to
be still interpreted as a {\it spatial probability density evolving in time
 with an appropriate similarity law}.

\vspace {0.25truecm}\noindent
\underbar{\it The Space-Time Fractional Diffusion Equation}

Recalling the approach by Mainardi, Luchko and Pagnini  in
\cite{Mainardi LUMAPA01}, to which we refer the interested reader for details,
it turns out that this
generalized diffusion equation, that we call
 {\it space-time fractional diffusion equation}, is
$$ {\, _x}D_{ \theta}^{\,\alpha} \,u(x,t) \, = \,
{\, _t}D_{*}^{\, \beta }\, u(x,t)\,,\q
-\infty<x<+\infty\,, \q t\ge 0\,,\eqno(4.2) $$
where  the $\alpha \,,\,\theta\,,\, \beta $ are real parameters
 restricted as follows
$$ 0<\alpha\le 2\,,\q |\theta| \le \min \{\alpha, 2-\alpha\}\,,\eqno(4.3)$$
$$ 0<\beta \le 1  \q\hbox{or} \q 1<\beta\leq \alpha \leq 2\,.\eqno(4.4) $$
Here
 $ {\,_x}D_{\,\theta}^{\,\alpha}$  and
${ \,_t}D_*^{\,\beta}  $  are  
integro-differential operators,
the {\it  Riesz-Feller space fractional derivative}
of order $\alpha $ and
asymmetry $\theta$
and the
{\it  Caputo time fractional derivative}
of order $\beta $, respectively.

The relevant  cases of the space-time fractional
diffusion equation  (4.2)
 include, in addition to the standard case of
{\it normal diffusion}   $\{\alpha =2, \beta =1\}$,
 the limiting case of the {\it D'Alembert wave} equation
 $\{\alpha =2, \beta =2\}$,
the {\it space fractional diffusion}
$\{0< \alpha <2,  \beta=1 \}$,
the {\it time  fractional diffusion}
$\{\alpha =2, 0< \beta <2\}$, and
the {\it neutral fractional diffusion}
$\{0<\alpha = \beta<2 \}$.
When $1<\beta \le 2$ we speak more properly about
the {\it fractional diffusion-wave equation} in that the corresponding
equation governs  intermediate phenomena between diffusion and wave
processes.

Let us now  resume the essential definitions
of the fractional derivatives in (4.2)-(4.4)
based on their Fourier and Laplace representations.

By denoting    the  Fourier  transform  of a  sufficiently well-behaved
(generalized) function  $f(x)$ as
$  \widehat f(\kappa)  =
{\cal F} \l\{ f(x);\kappa \r\}
  = \int_{-\infty}^{+\infty} \e^{\,\ds +i\kappa x}\,f(x)\, dx\,,
  \; \kappa \in \RR\,,$
 the {\it Riesz-Feller} space fractional derivative
of order $\alpha  $ and skewness $\theta$ is defined
as
$$\qq \qq \qq \qq \qq {\cal F} \l\{\, _xD_\theta^\alpha\, f(x);\kappa \r\} =
  - \psi_\alpha ^\theta(\kappa ) \,
  \, \widehat f(\kappa) \,,\qq \qq \qq \qq \qq \eqno(4.5)$$
 $$
   \psi_\alpha ^\theta(\kappa ) =
|\kappa|^\alpha \, \e^{\ds  i (\sgn \kappa)\theta\pi/2}
\,,\q
 0<\alpha  \le 2\,, \q
 |\theta| \le  \,\hbox{min}\, \{\alpha ,2-\alpha \}\,.
$$
 In other words the symbol of the pseudo-differential operator
$\,_xD_\theta^\alpha$ is  the logarithm of the
characteristic function of the generic strictly stable
probability density according to the Feller parameterization
\cite{Feller 52,Feller 71},
as  revisited by Gorenflo and Mainardi \cite{GorMai FCAA98}.
For this density we write
$$
            L_\alpha ^\theta (x) \Ft
\widehat{L_\alpha ^\theta} (\kappa ) =
     \exp \l[-\Psi_\alpha^\theta(\kappa )\r] \,,
\eqno(4.6) $$
where $\alpha $ is just the {\it stability exponent} ($0<\alpha \le 2$)
and $\theta$ is a real parameter related to the asymmetry
($ |\theta| \le  \,\hbox{min}\, \{\alpha ,2-\alpha \}$)
improperly called {\it skewness}.

By denoting   the Laplace transform of   a  sufficiently well-behaved
(generalized) function $f(t)$ as
$ \widetilde f(s) =
{\cal L} \l\{ f(t);s\r\}
 = \int_0^{\infty} \e^{\ds \, -st}\, f(t)\, dt\,, \;
\Re\,(s) > a_f\,,$
 the {\it Caputo} time fractional derivative  of
order $\beta  $ ($m-1 <\beta \le m\,,$ $\, m\in \NN$)
is defined through
\footnote{
The reader should observe that the {\it Caputo} fractional derivative
introduced in \cite{Caputo 67,Caputo 69,CaputoMaina 71}
represents a sort of regularization in the time origin for the
classical {\it Riemann-Liouville} fractional derivative
see \eg \cite{GorMai CISM97,Podlubny 99}.
We note that the {\it Caputo} fractional  derivative coincides with that
introduced  (independently of Caputo)
by  {\it Djrbashian \& Nersesian} \cite{DjrbashianNersesian 68},
 which has been   adopted by  Kochubei,
 see \eg \cite{Kochubei 89,Kochubei 90}
for treating initial value problems
in the presence of fractional derivatives.
In \cite{ButzerWestphal 00} the authors have pointed out
that such  derivative
was also considered by Liouville himself, but it should be noted that it
was disregarded by Liouville who did not recognize its role.}
$$ {\cal L} \l\{ _tD_*^\beta \,f(t) ;s\r\} =
      s^\beta \,  \widetilde f(s)
   -\sum_{k=0}^{m-1}    s^{\beta  -1-k}\, f^{(k)}(0^+) \,,
  \q m-1<\beta  \le m \,. \eqno(4.7)$$
This leads to define, see \eg \cite{GorMai CISM97,Podlubny 99},
$$
    _tD_*^\beta \,f(t) :=
\cases{
    {\ds \rec{\Gamma(m-\beta )}}\,{\ds\int_0^t
 {\ds {f^{(m)}(\tau)\, d\tau \over (t-\tau )^{\beta  +1-m}}}},
  & $ m-1<\beta  <m, $\cr\cr
     {\ds {d^m\over dt^m}} f(t)\,,
    & $\, \beta  =m\,. $\cr\cr }
   \eqno(4.8) $$

\vspace {0.25truecm}\noindent
\underbar{\it The Green Function}

When the diffusion equations (4.1), (4.2)  are equipped by
the initial and boundary conditions
$$u(x,0^+)=   \varphi(x)\,,
\q u(\pm\infty,t)=   0\,,  \eqno(4.9)  $$
their  solution reads
$ u(x,t) =
\int_{-\infty}^{+\infty} G(\xi ,t)\,
 \varphi (x-\xi) \, d\xi\,, $
where
 $G(x,t)$ denotes the fundamental solution (known as the {\it Green function})
corresponding to  $ \varphi (x) = \delta (x)$, the Dirac generalized
function.
We note that when $1<\beta \le 2$  we
must add a second initial condition $u_t(x, 0^+) = \psi(x)$,
which implies two Green functions corresponding to
$\{u(x,0^+) = \delta (x)$,   $u_t(x,0^+) = 0\}$
and
$\{u(x,0^+) = 0$, $u_t(x,0^+) = \delta (x)\}$.
Here we restrict ourselves
to consider  the first Green function because
only for this it is legitimate
to demand it to be a spatial probability density
evolving in time, see below.

It is straightforward to derive  from (4.2) the
composite Fourier-Laplace transform of the Green function by taking
into account the Fourier transform for the {\it Riesz-Feller}
space fractional derivative, see (4.5),,
and the Laplace transform for the {\it Caputo}
time fractional derivative,  see  (4.7).
We have, see \cite{Mainardi LUMAPA01} 
$$
  \Gks = {s^{\beta -1} \over
 s^\beta + \psi_\alpha^\theta (\kappa)} \,. \eqno(4.10)$$
By using the known scaling rules
for the Fourier  and Laplace  transforms,
we  infer 
without inverting the two transforms,
$$ \G(x,t)  =
    t^{-\gamma}\,\K (x/t^\gamma)\,,\q \gamma =\beta /\alpha
\,,   \eqno(4.11)     $$
where the one-variable function $\K$ is the
{\it reduced Green function} 
and $x/t^\gamma$ 
is the {\it similarity variable}.
We note  from
$$ \Gos = 1/s \lra \Got =1\,,$$
the {\it normalization property}
$$ \int_{-\infty}^{+\infty}\G(x,t) \,dx  =
     \int_{-\infty}^{+\infty}\K(x) \,dx =1\,,$$
and,  from
$$ \psi_\alpha^\theta (\kappa)
  = \overline{ \psi_\alpha^\theta (-\kappa)}
 =  \psi_\alpha^{-\theta} (-\kappa)
  \,,\,  $$
the {\it symmetry relation}
$$ K_{\alpha ,\beta}^\theta(-x) = K_{\alpha ,\beta}^{-\theta}(x)
\,,$$
 allowing us to restrict our attention to $x\ge 0\,.$

For $1<\beta \le 2$ we can show, see  \eg \cite{MainardiPagnini MELFI01},
that the second Green function is
a primitive (with respect to the variable $t$)
of the first Green function (4.11), so that, being no longer
normalized in $\RR$, cannot be interpreted as a spatial
probability density.

When $\alpha =2$ ($\theta=0$) and $\beta =1$
the inversion of the Fourier-Laplace
transform in (4.10) is trivial:
we recover the Gaussian density,
evolving in time with variance $\sigma^2 = 2t$,
typical of the normal diffusion,
$$ G_{ 2 ,1} ^{ 0}(x,t)
  =  {1 \over 2 \sqrt{\pi t}}\,
\exp\,\l( - x^2 / (4t)\r)    \,,\q x\in \RR\,,\q t >0\,,\eqno(4.12)$$
which exhibits the similarity law (4.11) with $\gamma =1/2$.
\newpage
For the analytical and computational determination of
the  reduced Green function
we  restrict our attention to $x>0\,$
because of the
{\it symmetry relation}. 
In this range Mainardi, Luchko and Pagnini \cite{Mainardi LUMAPA01}
have provided the Mellin-Barnes\footnote{
The names refer to the two authors, who in the first 1910's
developed the theory of these integrals  using them
for a complete integration of the hypergeometric differential equation.
However,
these integrals were first used by S. Pincherle in 1888.
For a revisited analysis of the pioneering work
by Pincherle (1853-1936, Professor of Mathematics at the
University of Bologna from 1880 to 1928) we refer to
Mainardi and  Pagnini \cite{MainardiPagnini OPSFA01}.}
integral representation
$$  \K(x) =
{1\over  \alpha x}
{1\over 2\pi i} \int_{\gamma-i\infty}^{\gamma+i\infty}
{\Gamma({s\over \alpha}) \, \Gamma(1-{s\over \alpha}) \,\Gamma(1-s)
 \over \Gamma(1-{\beta\over \alpha}s) \,
 \Gamma ( \rho \,s)\,
 \Gamma (1-\rho \,s)}
 \, x^{\,\ds s}\,  ds
\,,\; \rho =   { \alpha -\theta \over 2\,\alpha }\,,
\eqno(4.13)  $$
where $\gamma $ is a suitable real constant.
\vskip 0.25truecm 
\noindent
\underbar{{\it The Space Fractional Diffusion}} :
  $\{0<\alpha <2\,, \, \beta =1\}$

In this case we recover
the class $L_\alpha ^\theta(x)$ of the strictly stable (non-Gaussian)
densities exhibiting fat tails (with the algebraic decay proportional to
$ |x|^{-(\alpha +1)}$)
  and infinite variance,
$$  K_{\alpha ,1}^\theta (x) =   L_\alpha ^\theta(x) =
{1\over  \alpha x}
{1\over 2\pi i} \int_{\gamma-i\infty}^{\gamma+i\infty}
{\Gamma({s/ \alpha})  \,\Gamma(1-s)
 \over  \Gamma ( \rho \,s)\, \Gamma (1-\rho s)}
 \, x^{\,\ds s}\,  ds\,, \; \rho =   { \alpha -\theta \over 2\,\alpha }\,,
\eqno(4.14) $$
where
$0 <\gamma< \hbox{min} \{\alpha, 1 \}\,.$

A stable $pdf$ with extremal value of
the skewness parameter  is called {\it extremal}.
One can prove that all the extremal stable $pdf$s' with
$0<\alpha  <1$ are one-sided, the support being $\RR_0^+$  if
$\theta  =-\alpha \,, $
and  $\RR_0^-$  if $\theta  =+\alpha \,. $
The one-sided stable $pdf$'s with support in $\RR_0^+$ can be
better characterized by their (spatial) Laplace transform, which
turn out to be
$$ \widetilde {L_\alpha^{-\alpha}} (s) :=
   \int_0^\infty \!\! \e^{\ds\, -sx}\,L_\alpha^{-\alpha} (x) \,dx
 = \e^{\ds - s^\alpha} \,, \q \Re\, (s) >0\,, \q 0<\alpha <1\,.
\eqno(4.15)$$
 In terms of Mellin-Barnes integral representation we have
$$
  L_\alpha ^{-\alpha }(x) =
\frac{1}{ \alpha x}\,
\frac{1}{2\pi i} \, \int_{\gamma-i\infty}^{\gamma+i\infty}
 \frac{\Gamma({s/ \alpha})   }{ \Gamma ( s)}
 \, x^{\,\ds s}\,  ds\,, \q
0 <\gamma< \alpha< 1\,.
\eqno(4.16)$$

\vskip 0.25truecm \noindent
\underbar{{\it The Time Fractional Diffusion}} :
  $\{\alpha =2\,, \, 0<\beta <2\}$

In this case  we recover
the class
of the Wright type densities exhibiting stretched exponential
tails    and finite variance proportional to $t^\beta \,,$
$$  K_{2 ,\beta }^0 (x) = {1\over 2} M_{\beta /2}(x) =
{1\over  2 x}
{1\over 2\pi i} \int_{\gamma-i\infty}^{\gamma+i\infty}
{\Gamma(1-s)
 \over   \Gamma (1- \beta s/2)}
 \, x^{\,\ds s}\,  ds\,, \eqno(4.17) $$
where $0 <\gamma< 1\,.$
As a matter of fact ${1\over 2} M_{\beta /2}(x)$ turns out
to be a symmetric   probability density  related
  to the transcendental function  $M_\nu  (z)$
 defined for any  $\nu   \in (0,1)$ and $\forall z \in \CC$
as
$$ M_\nu  (z) =
  \sum_{n=0}^{\infty} 
  {(-z)^n\over n!\,\Gamma[-\nu   n + (1-\nu )]} =
  \rec{\pi} \sum_{n=1}^{\infty}{(-z)^{n-1} \over (n-1)!}
  \Gamma(\nu   n) \sin (\pi \nu   n).
  \eqno(4.18)    $$
We note that $M_\nu (z)$ is an entire function of
order $\rho =1/(1-\nu  )\,,$
that  turns out to be  a special case of the
Wright function\footnote{
 The Wright function   is defined by the  series  representation,
 valid in the whole complex plane,
  $$ \Phi _{\lambda ,\mu }(z ) :=
   \sum_{n=0}^{\infty}{z^n\over n!\, \Gamma(\lambda  n + \mu )}\,,
 \q \lambda  >-1\,, \,\q \mu \in \CC\,, \q z \in \CC\,.$$
Then, $ M_\nu (z) :=  \Phi _{-\nu , 1-\nu }(-z)$ with $ 0<\nu <1\,.$
The function $M_\nu(z)$ provides a generalization
of the Gaussian and of the Airy function in that
$$ M_{1/2}(z) = \rec{\sqrt{\pi}}\, \exp \l(-{\,z^2/ 4}\r)\,,\q
  M_{1/3}(z) =3^{2/3} \, {\hbox {Ai}} \l( {z/ 3^{1/3}}\r) \,.
$$}.
 Restricting   our attention to the
positive real axis ($r \ge 0$) we have:
the {\it Laplace transform}
  $$  \L \{M_\nu (r); s\} =
  E_\nu(- s )\,,\q \Re \,(s) >0\,, \q 0<\nu <1\,,\eqno(4.19) $$
 where $E_\nu $ is the Mittag-Leffler function,
and  the {\it asymptotic  representation} 
$$ \q\q\q\q\q\q\q\q M_\nu (r)  \sim  A_0  \,
     Y^{\, \nu - 1/2 } \, \exp \,\l( - Y\r)  \,,
\q x \to \infty\,, \q\q\q\q\q\q\q \eqno(4.20)$$
$$A_0 = {1\over   \sqrt{2\pi}\,(1-\nu )^\nu \, \nu ^{2\nu-1}} \,,\q
Y = (1-\nu )\, (\nu^\nu\, r)^{1/(1-\nu )}\,, $$
   a result formerly obtained by
Wright  himself,
and,   independently, by Mainardi and Tomirotti
\cite{MainardiTomirotti TMSF94} by using the saddle point method.
Because of the above exponential decay,  any moment
of order $\delta >-1$   for $M_\nu (r)\, $ is finite and given as
$$\int_0^{\infty}  r^{\,\ds \delta}\, M_\nu(r)\, dr
  = {\Gamma(\delta +1)\over \Gamma(\nu  \delta  +1)}\,,
   \q \delta > -1\,. \eqno(4.21)$$
In particular we get the normalization property in $\RR^+$,
$\int_0^\infty M_\nu (r)\, dr =1$.
In view of Eqs (4.11) and (4.21)
the moments (of even order)
of the fundamental solution
$G_{2,\beta}^0(x,t)$ turn out to be, for $n=0,1,2,\dots$
and $t \ge 0\,,$
$$ 
\int_{-\infty}^{+\infty} \!\!\! x^{2n}\,
   G_{2,\beta}^0(x,t)  \,dx =
  t^{\beta n}\, \int_0^\infty \!\!\! x^{2n} M_{\beta/2}(x)\, dx =
 {\Gamma(2n+1)\over \Gamma(\beta   n+1)} \, t^{\beta n}.
\eqno(4.22)$$

We agree to call $M_\nu (r)$ ($\nu  \in (0,1),\, r \in \RR_0^+$)
the {\it M-Wright function of order $\nu $},
understanding that its half represents  the spatial $pdf$
corresponding to the {\it time fractional diffusion} equation
of order $2 \nu$.
Relevant properties of this function,
see \eg \cite{GoLuMa 99,Mainardi CISM97,Mainardi LUMAPA01},
are concerning the limit expression for $\beta=1$, i.e.
$M_1(r) = \delta(r-1)$
and its relation with the extremal stable densities, i.e.
$$ {1 \over c^{1/\nu}}\,L_\nu^{-\nu }\l( {r\over c^{1/\nu}}\r)
 = {c\,\nu \over r^{\nu+1}}\,
    M_\nu\l( {c\over r^\nu }\r) \,,\q
  0<\nu  <1 \,,\q c>0\,,\q r>0
 \,. \eqno(4.23) $$

We note that, in  both  limiting cases of  space fractional  ($\alpha =2$)
and time fractional ($\beta =1$) diffusion, we recover the Gaussian density
of the {\it normal diffusion},
for which
$$
\begin{array}{ll}
 {\ds K_{2,1}^0(x)} & =  \, {\ds \frac{1}{  2 x} \,
\frac{1}{ 2\pi i}\, \int_{\gamma-i\infty}^{\gamma+i\infty}
 \frac{\Gamma(1-s)   }{\Gamma (1-s/2)}
 \, x^{\,\ds s}\,  ds}   
\\ \\
 & = \,  {\ds  L_{2}^0(x)} = 
 {\ds {1\over 2} M_{1/2} (x)
 = \frac{1}{2 \sqrt{\pi}}\, \e ^{\,\ds -x^2/4}}
 \,.
 \end{array}
 \eqno(4.24)
 $$
\vskip 0.25truecm \noindent
\underbar{{\it The Neutral Fractional Diffusion}} :
  $\{0<\alpha =\beta <2\}$

In this case,
surprisingly,
the corresponding (reduced)
Green function     can be expressed (in explicit form)
in terms of a  (non-negative) simple elementary function,
that we denote  by $   N_{ \alpha} ^{ \theta}(x)\,,$
as it is shown in  \cite{Mainardi LUMAPA01}:
$$  
\begin{array}{ll}
{\ds N_{ \alpha} ^{ \theta}(x)} & =
\, {\ds {1\over  \alpha x}\,
{1\over 2\pi i} \, \int_{\gamma-i\infty}^{\gamma+i\infty}
{\Gamma({s\over \alpha}) \, \Gamma(1-{s\over \alpha})
 \over \Gamma ( \rho \,s)\, \Gamma (1-\rho \,s)}
 \, x^{\,\ds s}\,  ds} \\ \\
& =\, {\ds {{1\over\pi}}\,{x^{\alpha-1} \sin[{\pi\over 2}(\alpha -\theta )] \over
1 + 2x^\alpha \cos[{\pi\over 2}(\alpha -\theta)] + x^{2\alpha}}}
 \,.
 \end{array}
 \eqno (4.25)
 $$

\vspace {0.25truecm}\noindent
\underbar{\it The Fractional Diffusion  Processes}

The {\it self-similar stochastic processes}
generated by the above probability densities evolving in time
can be considered as generalizations
of the standard diffusion processes and therefore
distinguished from it with the label "fractional".
When $0<\beta <1  $  random walk models
can be  introduced
to generalize the classical  Brownian motion of the standard diffusion,
as it was investigated  in a number of papers of our group
see \eg  \cite{GorMai FCAA98,GorMai ChemPhys02}.
In the case of {\it space fractional diffusion}
we obtain a special class of Markovian processes, called
stable L\'evy motions, which exhibit infinite variance
associated to the possibility of arbitrarily large jumps ({\it L\'evy
flights}).  In the  case of {\it time fractional diffusion}
 we obtain a  class of stochastic processes which
 are non-Markovian
and exhibit a variance consistent with slow anomalous diffusion.
For the general genuine {\it space-time fractional} diffusion 
  ($0<\alpha <2\,,\,0<\beta  < 1$)
we generate   a class of  densities
(symmetric or non-symmetric according to $\theta =0$ or $\theta \ne 0$)
which exhibit  fat tails
with an algebraic decay $\propto |x|^{-(\alpha +1)}\,.$
Thus they belong to the domain of attraction of the L\'evy stable densities
of index $\alpha $ and  can be referred to as
{\it fractional stable densities}.
The related stochastic processes
possess the  characteristics of the previous two classes;
indeed,  they are non-Markovian (being $0< \beta < 1$)
and exhibit infinite variance
associated to the possibility of arbitrarily large jumps
(being $0< \alpha <2$).

\vskip 0.5truecm

\section*{5. Subordination for space fractional diffusion processes}

In the book by Feller, see \cite{Feller 71}, at p. 176 we read:
{\it Let $X$ and $Y$ be independent strictly stable variables,
with characteristic exponent $\alpha$ and $\beta $ respectively.
Assume $Y$ to be a positive variable (whence $\beta <1$).
The product $X\, Y^{1/\alpha }$   has a stable distribution
with exponent $\alpha \, \beta $.}

In other words, 
this statement means that any strictly stable
process (of  exponent $\gamma =\alpha \cdot \beta$)
is subordinated to a parent strictly stable process (of exponent $\alpha $)
and directed by an extremal  
strictly stable process (of exponent $0<\beta <1$).
Feller's  proof is vague being, as a matter of fact,
 limited to {\it symmetric} subordinated and parent
stable     distributions. Furthermore, the proof,   scattered in
several sections, is essentially based on the use
of Fourier and Laplace transforms.    Here we would like to
make more precise the previous statement by Feller by considering
the possibility of asymmetry  characterized by the index $\theta$
as previously explained and making use of
the Mellin machinery outlined in Sections 2 and 3,
and  of the
Mellin-Barnes integral representation (4.14).
So, in virtue of the Mellin inversion formula (2.3)
we can write for the generic {\it strictly stable} $pdf$
$$
L_{\alpha}^{\theta}(x) \Mt \frac{1}{\alpha}
\frac{\Gamma\lt(-\frac{s-1}{\alpha} \rt)\Gamma[1 + (s-1)]}
{\Gamma[1+\rho (s-1)]\Gamma[-\rho(s-1)]}\,,
\q \rho =   { \alpha -\theta \over 2\,\alpha }\,.
\eqno(5.1) $$
Let us now consider the evolution in time according
to the {\it space fractional diffusion} equation by writing
$$     L_{\alpha}^{\theta}(x; t):= G_{\alpha,1}^\theta (x,t)
  = t^{-1/\alpha }\,
          L_{\alpha}^{\theta}\l( \frac{x}{t^{1/\alpha}}\r)\,. \eqno(5.2)$$
Then  we prove the following 

\noindent {\sc Theorem}.
\it
Let $L_{\alpha_p}^{\theta_p}(x;t),$ $L_{\alpha_q}^{\theta_q}(x;t)$
and $L_{\beta}^{\theta_{\beta}}(x;t)$ be strictly stable densities
with exponents
$\alpha_p \,, \alpha_q \,, \beta$ and asymmetry parameters
$\theta_p \,, \theta_q \,, \theta_{\beta}$, respectively,
such that
$$ 0 < \alpha_p \leq 2 \,,  \q
|\theta_p| \leq \min\{\alpha_p, 2-\alpha_p\} \,,$$
$$ 0 < \alpha_q \leq 2 \,,    \q
|\theta_q| \leq \min\{\alpha_q, 2-\alpha_q\}\,,$$
$$ 0<\beta \leq 1 \,,    \q
\theta_{\beta} = -\beta \,,$$
then the following subordination formula 
holds true for $0 < x < \infty\,,$
$$
L_{\alpha_p}^{\theta_p}(x;t)=\int_0^{\infty}
L_{\alpha_q}^{\theta_q}(x;\tau)\, L_{\beta}^{-\beta}(\tau;t)
\, d\tau\,, \; \hbox{with} \;\;
\alpha_p=\beta \alpha_q \,,\; \theta_p=\beta \theta_q\,.
\eqno(5.3) $$
\rm
Because of the scaling property (5.2) of the stable $pdf$'s,
we can alternatively state
$$
t^{-1/\alpha_p}L_{\alpha_p}^{\theta_p}
\lt(\frac{x}{t^{1/\alpha_p}}\rt)=\int_0^{\infty}
\tau^{-1/\alpha_q}L_{\alpha_q}^{\theta_q}\lt(\frac{x}{\tau^{1/\alpha_q}}\rt)
t^{-1/\beta}L_{\beta}^{-\beta}\lt(\frac{\tau}{t^{1/\beta}}\rt)
\, d\tau \, .
\eqno(5.4)  $$


The proof of Eq. (5.4) is a (straightforward)
consequence of the previous considerations.
By recalling the Mellin pairs for the involved stable densities
(that can be easily obtained from Eq. (5.1)-(5.2)
 by adopting the correct parameters) 
and the scaling properties of the Mellin transform, see (2.4),
we have
$$
t^{-1/\alpha_p} L_{\alpha_p}^{\theta_p}
\lt(\frac{x}{t^{1/\alpha_p}}\rt) \Mt
t^{-1/\alpha_p}\lt(\frac{1}{t^{1/\alpha_p}}\rt)^{-s}
\frac{1}{\alpha_p}
\frac{\Gamma\lt(-\frac{s-1}{\alpha_p} \rt)\Gamma[1 + (s-1)]}
{\Gamma[1+\rho_p(s-1)]\Gamma[-\rho_p(s-1)]}
\eqno(5.5) $$
and
$$
b \, c\, x^a \, L_{\beta}^{\theta_{\beta}}(c x^b) \Mt
 c^{1- \frac{s+a}{b}} \frac{1}{\beta}
\frac{\Gamma\lt(-\frac{s+a}{b \beta}+\frac{1}{\beta} \rt)
\Gamma\lt[1 + \lt(\frac{s+a}{b}-1\rt)\rt]}
{\Gamma\lt[1+\rho_{\beta}\lt(\frac{s+a}{b}-1\rt)\rt]
\Gamma\lt[-\rho_{\beta}\lt(\frac{s+a}{b}-1\rt)\rt]} \,.
\eqno(5.6)$$
After some algebra  we recognize
$$
\M\{t^{-1/\alpha_p}\, L_{\alpha_p}^{\theta_p}
\lt(\frac{x}{t^{1/\alpha_p}}\rt);s\}=
\M\{b x^a \, c \, L_{\beta}^{\theta_{\beta}}(c x^b);s\}
\, \M\{L_{\alpha_q}^{\theta_q}(x);s\}
\eqno(5.7) $$
provided that
$$
\theta_{\beta}=-\beta \,,
\q a=\alpha_q-1\,, \q b=\alpha_q \,,\q c=t^{-1/\beta}\,,
\eqno(5.8)
$$
and
$$
\alpha_p=\beta \alpha_q \,,\q \theta_p=\beta \theta_q \,.
\eqno(5.9) $$
Recalling the Mellin convolution formula  (2.5)
we obtain from Eqs. (5.7)-(5.9)
the integral representation
$$
t^{-1/\alpha_p}\, L_{\alpha_p}^{\theta_p}
\lt(\frac{x}{t^{1/\alpha_p}}\rt)=\int_0^{\infty}
\alpha_q \xi^{\alpha_q-1}\, L_{\alpha_q}^{\theta_q}\lt(\frac{x}{\xi}\rt)
t^{-1/\beta}\, L_{\beta}^{-\beta}\lt(\frac{\xi^{\alpha_q}}{t^{1/\beta}}\rt)
\, \frac{d\xi}{\xi} \, ,
\eqno(5.10)
$$
and, by replacing  $\xi  \rightarrow \tau^{1/\alpha_q}$,
we finally get Eq. (5.4).

Taking into account the relationships in  Eq. (5.9) 
we can point out  some interesting {\it subordination laws}.
In particular we  observe that any {\it symmetric} stable distribution
with exponent $\alpha_p =\alpha \in (0,2)$
($\theta_p =\theta =0$) is subordinated to
the {\it Gaussian distribution} ($\alpha_q=2\,,\, \theta_q=0$),
see $L_2^0(x)$ in (4.24), through an extremal stable density of exponent
$\beta =\alpha /2$, that is
$$
L_{ \alpha }^{0}(x;t)= \int_0^\infty
L_{2}^{0}(x;\tau)\, L_{\alpha /2}^{-\alpha/2 }(\tau;t)\, d\tau \,,
\q 0<\alpha <2\,.
\eqno(5.11)$$
Furthermore, by recalling
the {\it generalized Cauchy density}
of skewness $\theta$ ($|\theta|<1$),
see \eg Eq. (4.9) in \cite{Mainardi LUMAPA01},
$$
 L_1^{\theta} (x) =
 \frac{1}{\pi} \, \frac{  \cos (\theta \pi/2)
 }{ [x+   \sin (\theta \pi/2)]^2 +[ \cos (\theta \pi/2)]^2 }\,,
\; |\theta| <1\,, \; -\infty < x <+\infty\,,
\eqno(5.12) $$
we note that
 any stable distribution
with exponent $\alpha_p = \alpha \in (0,1)$
and skewness $|\theta_p|= |\theta| <\alpha $
is subordinated to the {\it generalized Cauchy
distribution} with skewness $|\theta_q| = |\theta|/\alpha  <1$
denoted by  $L_1^{\theta/\alpha}(x)$, see (5.12),
through an extremal density
of exponent $\beta =\alpha $, that is
$$ L_{\alpha}^{\theta}(x;t)= \int_0^\infty
L_{1}^{\theta/\alpha }(x;\tau)\,  L_{\alpha }^{-\alpha }(\tau;t)\, d\tau \,,
\q 0<\alpha <1\,.
\eqno(5.13)$$

\section*{6. Subordination for time fractional diffusion processes}

For the {\it M-Wright function} we deduce from (4.17)
 the Mellin transform pair:
$$
     M_{\nu }(r) \Mt
\frac{\Gamma\lt[1 +(s-1)\rt]}{\Gamma[1 + \nu (s-1)]}  \,,
\q 0<\nu <1\,.\eqno(6.1)$$
Let us now consider the evolution in time according
to the corresponding {\it time  fractional diffusion} equation by writing
$$     M_{\nu }(x; t) := 2\, G_{2,2\nu}^0 (x,t)
= t^{-\nu  }\,  M_{\nu }\l( \frac{x}{t^{\nu }}\r)\,. \eqno(6.2)$$
Then  we prove the following 

\noindent {\sc Theorem}.
\it
Let $M_{\nu }(x;t),$ $M_{\eta}(x;t)$
and $M_{\beta}(x;t)$ be $M$-Wright functions
of  orders
$\nu, \eta, \beta \in (0,1)$
 respectively,
then the following subordination formula 
holds true for $0 < x < \infty$,
$$
M_{\nu}(x,t) =\int_0^{\infty}
M_{\eta }(x;\tau)\, M_{\beta}(\tau;t)
\, d\tau\,, \; \hbox{with} \;\;
\nu = \eta \, \beta\,.
\eqno(6.3) $$
\rm
Because of the scaling property (6.2) of the $M$-Wright functions,
we can alternatively state
$$
t^{-\nu }\, M_{\nu }
\lt(\frac{x}{t^{\nu }}\rt)=\int_0^{\infty}
\tau^{-\eta }\, M_{\eta }\lt(\frac{x}{\tau^{\eta}}\rt)
t^{-\beta}\, M_{\beta} \lt(\frac{\tau}{t^{\beta}}\rt)
\, d\tau \, .
\eqno(6.4)  $$
The proof of Eq. (6.4) is a (straightforward) consequence
of the previous considerations.
After some algebra  we recognize
$$
\M\{t^{-\nu }\, M_\nu
\lt(\frac{x}{t^{\nu }}\rt);s\}=
\M\{b\,c\,  x^a \, M_{\beta}(c\, x^b);s\}
\, \M\{M_{\eta}(x);s\}
\eqno(6.5) $$
provided that
$$
 a=\frac{1}{\eta}-1 \,, \q b=\frac{1}{\eta}\,,\q c=t^{-\nu /\eta}\,,
\eqno(6.6)
$$
and
$$
   \nu = \beta \,  \eta \,.
\eqno(6.7) $$
Recalling the Mellin convolution formula  (2.5)
we obtain from Eqs. (6.5)-(6.7)
the integral representation
$$
t^{-\nu }\, M_{\nu }
\lt(\frac{x}{t^{\nu }}\rt)=\int_0^{\infty}
\frac{1}{\eta} \, \xi^{1/\eta -1}\, M_{\eta }\lt(\frac{x}{\xi}\rt)
t^{-\beta}\, M_{\beta}\lt(\frac{\xi^{1/\eta }}{t^{\beta}}\rt)
\, \frac{d\xi}{\xi} \, ,
\eqno(6.8)
$$
and, by replacing  $\xi  \rightarrow \tau^{\eta}$,
we finally get Eq. (6.4).

We note that for $\eta=1/2$  the corresponding $M$-Wright function
reduces to  twice  the Gaussian density according to (4.24), so the
subordination formula (6.3) reads
$$ M_{\beta/2} (x;t) = 2 \int_0^\infty L_2^0(x;\tau)\, M_\beta(\tau;t)
\, d\tau\,. \eqno (6.9)$$


\section*{7. Subordination for space-time fractional diffusion}

For the {\it space-time fractional} diffusion we shall
prove two relevant {\it subordination laws}.
For this purpose we first consider the evolution in time writing
$$     K_{\alpha,\beta }^{\theta}(x; t):= G_{\alpha,\beta }^\theta (x,t)
  = t^{-\beta /\alpha }\,
   K_{\alpha,\beta }^{\theta}\l( \frac{x}{t^{\beta/\alpha}}\r)\,. \eqno(7.1)$$
 Then the required subordination laws read
$$
K^{\theta}_{\alpha,\beta} (x;t) =
\int_0^{\infty} L_{\alpha}^{\theta} (x;\tau)
\, M_{\beta}(\tau;t)\,d\tau
\,, \q 0<\beta \le1 \,. \eqno(7.2)$$
and
$$
K^{\theta}_{\alpha,\beta} (x;t) =
\int_0^{\infty} N_{\alpha}^{\theta} (x;\tau)
\, M_{\beta/\alpha }(\tau;t)\,d\tau
\,, \q 0<\beta/\alpha \le 1 \,. \eqno(7.3)$$
In order to prove the above laws we start
from two relevant results in \cite{Mainardi LUMAPA01}, see (6.1),
that we report below
$$ \K(x) = \cases{
   \alpha   \, {\ds \int_0^\infty} \!\!
  \l[\xi ^{\alpha -1}\,  {M}_{\beta}\l(\xi ^{\alpha}\r)\r]\,
  {L}_{\alpha}^\theta\l({x/\xi }\r) \,
      {\ds{d\xi  \over \xi }} \,,&$\; 0<\beta \le 1\,,$ \cr\cr
  {\ds \int_0^{\infty}} \!\!
M_{\beta / \alpha}(\xi ) \,N_\alpha^\theta(x/\xi)\,
    {\ds{d\xi  \over \xi} }\,, & $\; 0<\beta/\alpha\le 1 \,.$ \cr}
\eqno(7.4)$$
In \cite{Mainardi LUMAPA01}
the above identities  have enabled us
to  {\it extend the probability interpretation
of the Green functions to the ranges}
$\{0<\alpha < 2\} \cap \{0<\beta < 1 \}\, $
{\it and} $\, \{1<\beta < \alpha < 2\}\,.$
Indeed, the formulae in (7.4) show
the non-negativity of the Green functions of the genuine space-time
fractional  diffusion equation based on the non-negativity of the
Green functions of the particular cases of space,  time and neutral
fractional diffusion.
We note that the  formulae were derived by using the Mellin-Barnes
representation
of the corresponding Green functions,
a method akin with that of the  Mellin transform  as we have
previously seen.
Let us now prove Eq. (7.2).
Because of the scaling properties of the involved functions
this is equivalent to  prove, with $0<\beta \le 1$,
 $$t^{-\beta /\alpha}\,
K^{\theta}_{\alpha,\beta} \l( {x \over t^{\beta /\alpha}}\r )
=\int_0^{\infty} \!\!
 \tau ^{-1/\alpha }\,
L_{\alpha}^{\theta} \l( {x \over \tau^{1/\alpha}} \r )  \,
  t^{-\beta }\, M_{\beta}\l ({\tau \over t^\beta} \r )\, d\tau \,,
\eqno(7.5)$$
that can be easily achieved 
from the first equation of (7.4),
by replacing $x$ with  $x/t^{\beta/\alpha}$
and making  the change of variable $\xi  = \tau^{1/\alpha} /t^{\beta/\alpha} $
in the integral.
Similarly, to prove Eq. (7.3) we can verify, 
with  $0<\frac{\beta}{\alpha} \le 1$,
$$t^{-\beta /\alpha}\,
 K^{\theta}_{\alpha,\beta} \l( {x \over t^{\beta /\alpha}}\r )
=\int_0^{\infty} \!\!
 \tau ^{-1}\,
 N_{\alpha}^{\theta} \l ({x\over \tau}\r )  \,
  t^{-\beta/\alpha }\,
 M_{\beta/\alpha } \l( {\tau \over t^{\beta/\alpha}} \r)\, d\tau ,
 \eqno(7.6)
$$
In this case it suffices to replace $x$ with  $x/t^{\beta/\alpha}$
and   set $\xi  = \tau /t^{\beta/\alpha}  $
in the second equation of (7.4).
\newpage

It is worth to note that whereas for the particular cases of space and time
fractional diffusion  the corresponding subordination
formulas (5.3) and (6.3) involve functions of the same class,
here the subordination formulas (7.2) and (7.3)
involve functions of different classes. However, in all the cases
we point out that {\it the $pdf$ of the directing process
is a M-Wright function}, including  (5.3)
if we take into account  (4.23).

If we like  to use the terminology  $\G(x;t)$ for  the  Green functions,
the two subordination laws
(7.2)-(7.3) can be resumed as follows
$$
G^{\theta}_{\alpha,\beta}(x;t)=
 \cases{
  2   \, {\ds \int_0^{\infty} G_{\alpha,1}^{\theta}(x;\tau)\,
  G_{2,2\beta}^0(\tau;t)\, d\tau }  \,,
     &$\; 0<\beta \le 1\,,$ \cr\cr
  2\, {\ds \int_0^{\infty} G_{\alpha,\alpha }^{\theta}(x;\tau)\,
 G_{2,2\beta/\alpha}^0(\tau;t)\, d\tau }\,,
      &$\; 0<\beta/\alpha  \le 1\,.$ \cr\cr }
                   \eqno(7.7)$$
From the first equation in (7.7)
we note that
the  Green function for
the space-time fractional diffusion equation of
order $\{\alpha,\beta \}\,,$
with $0<\alpha \le 2$ and $0<\beta \le 1\,,$
can be expressed in terms  of the  Green function for
the space fractional diffusion equation of order $\alpha$
and the  Green function for
the time fractional diffusion equation of order $2\beta\,. $

\section*{8. Concluding discussion}

There are several ways of deriving subordination formulas for
fractional diffusion processes.
A natural way is to
to start from an approximating continuous time random walk model and
carry out an appropriate passage to the limit. As we cannot give here
a comprehensive survey of achievements of contributors to this
subject, let us only cite \cite{BaeumerMeerschaert 01,Nevada PRE02}
 for a concise description.
In contrast to this "stochastic" way it was our aim to show how
classes of subordination formulas can be found by purely analytic
methods, namely by exploiting solution formulas for the fundamental
solutions of the essential types of fractional diffusion equations,
formulas  expressing these solutions in form of Mellin-Barnes
integrals (see \cite{Mainardi LUMAPA01}).
Using the machinery offered by Mellin transform
theory these formulas can be rewritten as integrals of products that
in turn allow a probabilistic interpretation as subordination
formulas. So, we testify for the fact that the fascinating field of
fractional diffusion processes  is not only interesting from the
view-points of probability theory and statistical physics, but also
from that of pure analysis where it is a playground for
lovers of special functions
and integral transforms.
In particular, we also demonstrate that the Mellin transform
is an important and useful integral transform in its own right.

\vfill\eject


\end{document}